\documentclass{article}
\usepackage[centertags]{amsmath}
\usepackage{amsfonts}
\usepackage{amssymb}
\usepackage{amsthm}
\usepackage{epsfig, amssymb, subfigure}

\usepackage[ left=1in, right=1in, top=1in, bottom=1in, includefoot]{geometry}

\begin{document}

\newtheorem{theorem}{Theorem}

\newtheorem{proposition}{Proposition}

\newtheorem{lemma}{Lemma}

\theoremstyle{definition}
\newtheorem{definition}{Definition}
\theoremstyle{remark}
\newtheorem{corollary}{Corollary}
\newtheorem{Remark}{Remark}
\newtheorem{Problem}{Problem}




\def\STU#1#2{      
\begin{array}{ccccc}
\epsfig{file=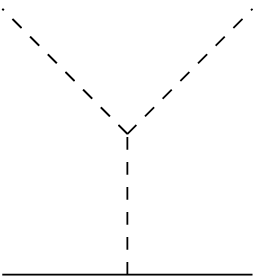, height=#1 } & \raisebox{#2}{=} &\epsfig{file=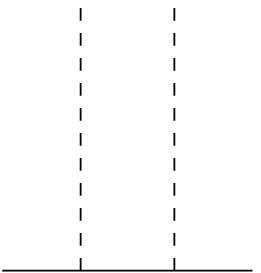, height=#1 } & \raisebox{#2}{-} & \epsfig{file=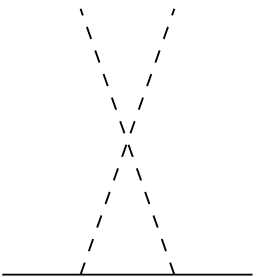, height=#1 }
\end{array}}

\def\IHX#1#2{      
\begin{array}{ccccc}
\epsfig{file=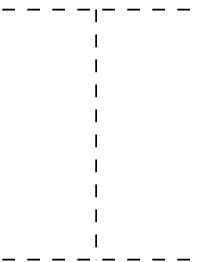, height=#1 } & \raisebox{#2}{=} &\epsfig{file=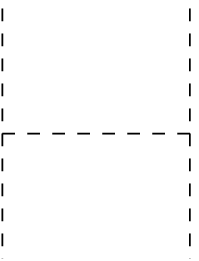, height=#1 } & \raisebox{#2}{-} & \epsfig{file=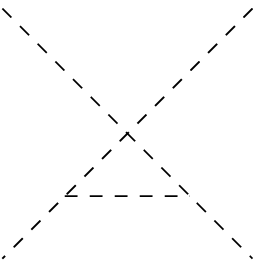, height=#1 }
\end{array}}

\def\AS#1#2{      
\begin{array}{ccc}
\epsfig{file=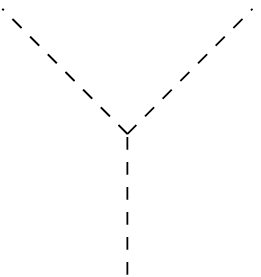, height=#1 } &  \raisebox{#2}{=} & \raisebox{#2}{-} \epsfig{file=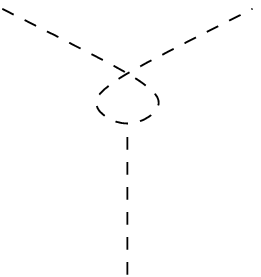, height=#1 }
\end{array}}


\def\bracx#1#2{\langle #1 ,  #2 \rangle_{X} }
\def\FG{\sideset{}{_{X_M}^{FG}}\int}
\def\fgi{\int_{X_M}^{FG}}
\def\bracxm#1#2{\langle #1 ,  #2 \rangle_{X_M} }
\def\bracxmt#1#2{\langle #1 ,  #2 \rangle_{X_M}^{t} }
\def\brat#1#2{\langle #1 ,  #2 \rangle^{t}  }
\def\bracy#1#2{\langle #1 ,  #2 \rangle_{Y} }
\def\bracs#1#2{\langle #1 ,  #2 \rangle_{S} }


\def\paamlt{Z_{0}^{M:t} (\sml )}
\def\aamlt{Z^{M:t} (\sml )}
\def\aart#1{Z^{M;t} (#1)}
\def\paart#1{Z_{0}^{M;t} (#1)}
\def\aar#1{Z^M (#1)}
\def\aaml{Z^M (\sml )}
\def\zm#1{Z^M (#1)}
\def\paaml{Z_{0}^{M} (\sml)}
\def\paar#1{Z_{0}^{M} (#1)}
\def\ea{Z^{M;h,t}(\lambda_i \otimes 1_M)}
\def\eb{\pi^{h,t} (Z^M (1_0 \otimes \sml )^{-1} D_i(Z^M ( \sml ))^{b_i})}
\def\Zmht{Z^{M;h,t}}
\def\zmht#1{Z^{M;h,t} (#1)}
\def\bra#1#2#3{\langle #1 ,  #2 \rangle_{#3} }
\def\pzm#1{Z_{0}^{M} (#1)}
\def\z{\check{Z}}
\def\zmhtn#1{Z_{\leq n}^{M;h,t} (#1)}
\def\Zhtn#1{Z_{\leq n}^{h,t} (#1)}
\def\zht#1{\check{Z}^{h,t} (#1)}
\def\pzm{Z^{M}_{0}}
\def\Zm{Z^{M}}
\def\zm#1{Z^{M} (#1)}
\def\zmt#1{Z^{M;t} (#1)}
\def\ZM{Z^{M}}
\def\Zht#1{Z^{h,t} (#1)}
\def\ZMH{Z^{M;h}}
\def\PZMT{Z_{0}^{M;t}}
\def\ZMT{Z^{M;t}}
\def\ZM{Z^{M}}
\def\PZM{Z^{M}_{0}}

\def\sm{\sigma_{M}}
\def\pa#1{\partial ( #1 )}
\def\bsig{B^{\Sigma}}
\def\h#1#2{H_{#1}(#2)}
\def\fun#1{\pi_1 (#1)}
\def\pres#1#2{\langle #1 | #2 \rangle}
\def\map#1#2#3{#1 : #2 \rightarrow #3}
\def\quo#1{#1 / #1_{q+1} }
\def\fquo#1{\frac{#1}{#1_{q+1}} }
\def\fquoq#1{\frac{#1}{#1_{q+1}} }
\def\art{\textit{Art}}
\def\pow#1{\mathcal{P}$$ \textit{(#1)} }
\def\ser#1#2{ #1_1 , \ldots , #1_{#2}}
\def\MU#1#2{\mu_{i}^{(n)}}
\def\cp{\vartriangle}
\def\Exp#1{\exp (#1)}
\def\coeff#1#2{\text{Coeff}(#1,#2)}
\def\prodab#1#2{\prod_{i= #1}^{#2} a_{i}^{b_i}}
\def\ab#1{a_{#1}^{b_{#1}}}
\def\lmi{\lambda_{i}^{M}}
\def\li{\lambda_{i}}
\def\Exp#1{\text{exp} (#1)}
\def\pht{\pi^{h,t}}
\def\sml{\sigma_{LM}}
\def\d2i{D^2 \times I}
\def\sxm{\sigma_{X_M}}
\def\mubar{\bar{\mu}}
\def\chiso{\raisebox{1mm}{$\chi$}}
\def\quo#1{#1 / #1_{n+1} }
\def\quoq#1{#1 / #1_{q+1} }
\def\bigbrac#1{\left(\begin{array}{c}#1 \end{array}\right)}


\def\Cht#1{\mathcal{C}$$^{h,t} (#1) }
\def\Ct#1{\mathcal{C}$$^{t} (#1) }

\def\Alm{\mathcal{A}$$(\uparrow_{X_L \cup X_M})}
\def\alglm{\mathcal{A}$$(\uparrow_{X_L }, X_M)}
\def\algle{\mathcal{A}$$(\uparrow_{X_L }, \emptyset )}
\def\Al1m{\mathcal{A}$$(\uparrow_{X_L +1 \cup X_M})}
\def\algl1e{\mathcal{A}$$(\uparrow_{X_L+1 }, \emptyset)}
\def\alghtl1e{\mathcal{A}$$^{h,t}(\uparrow_{X_L+1 }, \emptyset)}
\def\algyz{\mathcal{A}$$(\uparrow_X, Y)}
\def\algspecial{\mathcal{A}$$^{h,t}_{n-[ \frac{n}{2} ]-1}(\uparrow_{Y_0 \cup Y_L }, \emptyset)}
\def\Ahtlo{\mathcal{A}$$^{h,t}(\uparrow_{X_L +1 })}
\def\Bhtlo{\mathcal{B}$$^{h,t}(X_L +1 )}
\def\algtl1e{\mathcal{A}$$^t(\uparrow_{X_L+1 }, \emptyset)}
\def\Ahtl1m{\mathcal{A}$$^{h,t}(\uparrow_{X_L +1 \cup X_M})}
\def\Atl{\mathcal{A}$$^t(\uparrow_{X_L })}
\def\algtle{\mathcal{A}$$^{t}(\uparrow_{X_L }, \emptyset)}
\def\algez{\mathcal{A}$$(\emptyset, Y)}
\def\algyzn{\mathcal{A}$$_n(\uparrow_{X}, Y)}
\def\algchiimage{\mathcal{A}$$(\uparrow_{X \cup Y^{\prime}}, Y-Y^{\prime})}
\def\algyzh{\mathcal{A}$$^h(\uparrow_{X}, Y)}
\def\algyzt{\mathcal{A}$$^t(\uparrow_{X}, Y)}
\def\algrs{\mathcal{A}$$(\uparrow_X, Y)}
\def\algqq#1#2{\mathcal{A}$$(\uparrow_{#1 }, #2)}
\def\algxy{\mathcal{A}$$(\uparrow_X, Y)}
\def\algxyt{\mathcal{A}$$^t(\uparrow_X, Y)}
\def\algxyh{\mathcal{A}$$^h(\uparrow_X, Y)}
\def\Ahoxl1{ \mathcal{A}$$^{h(0)}(\uparrow_{X_L +1})}
\def\Axl1{ \mathcal{A}$$(\uparrow_{X_L +1})}
\def\algtlm{\mathcal{A}$$^{t}(\uparrow_{X_L }, X_M)}

\title{The \AA rhus Integral and the $\mu$-invariants}
\author{Iain Moffatt\footnote{
Mathematics Institute, University of Warwick, Coventry, CV4 7AL, UK.
\newline ${\hspace{.35cm}}$ {\em Current address:} Department of Applied Mathematics (KAM), Charles University, Prague, Czech Republic;
\newline ${\hspace{.35cm}}$ email: \texttt{iain@kam.mff.cuni.cz.}
\newline ${\hspace{.35cm}}$ This version: 4th February 2005. ${\hspace{0.75cm}}$ First version: 3rd March 2004
}}
\date{}
\maketitle

\abstract{We relate the tree part of the \AA rhus integral to Milnor's $\mu$-invariants of string-links in homology balls thus generalizing results of Habegger and Masbaum.}

\section{Introduction}

 Milnor's $\mubar$-invariants of Links and their well defined cousins, the $\mu$-invariants of string-links are classical and well studied invariants.
These invariants have been brought into the realm of finite-type invariants by Bar-Natan in \cite{BN:95}, Lin in \cite{Li:97} and Habegger and Masbaum in \cite{HM:00}.
In this paper we are particularly interested in Habegger and Masbaum's formula which expresses the $\mu$-invariants
 in terms of the tree part of the Kontsevich integral.

  The literature on the $\mubar$-invariants is mostly concerned with links in $S^3$.  The generalization to $\mubar$-invariants of links in integral homology spheres and  $\mu$-invariants of string-links in homology balls exists mostly as folklore.  We will discuss the $\mu$-invariants of string-links in homology balls and generalize Habegger and Masbaum's results by relating the $\mu$-invariants to the \AA rhus integral, which is a generalization of the Kontsevich integral to links in rational homology spheres, defined in \cite{BGRT:AI, BGRT:AII, BGRT:AIII}.
We do this by representing string-links in homology spheres by string-links in $\d2i$ with some distinguished surgery components.
Note that the \AA rhus integral agrees, up to a normalization, with the LMO invariant.

A familiarity with the basic properties of the Kontsevich integral and the relevant algebras is assumed.

\section{Tangles and String-links} \label{sec:defs}

Let $B^M$ be a connected, compact orientable 3-manifold equipped with a fixed identification $\varphi$ of  its boundary with $\pa{\d2i}$.
A {\em tangle} of $n$ components $T \subset B^M$ is a smooth compact 1-manifold $X$, of $n$ components, together with a smooth  embedding
$T: (X, \partial (X) ) \rightarrow (B^M , \partial (B^M))$, transverse to the boundary.
As is standard, we abuse notation and  confuse a tangle, its embedding and its isotopy class.

By a \textit{framing} on a component $i$  of a tangle we mean that we equip  $i$ with a non-vanishing vector field such that the restriction to the boundary is the restriction of a fixed unit vector field normal to the $x$-axis of $D^2$ under the identification.

A {\em coloured} tangle is a tangle equipped with a bijection from the set of components onto a set of cardinality $n$, where $n$ is the number of components.

Since $B^M$ can be obtained by surgery on a framed link $L \subset \d2i$, we may represent a tangle $T \subset B^M $ by a tangle $T^{\prime} \subset \d2i$ some of whose components are distinguished framed copies of $S^1$, on which we do the surgery. We say that $T^{\prime}$ {\em represents} $T$.
We will call these distinguished components the {\em surgery components} and the other components the {\em linking components}.
If the tangle $T$ is coloured then this partitions  the colouring set into  sets  corresponding to the surgery components and the linking components. In this paper we denote these  sets $X_M$ and $X_L$ respectively.

\begin{figure}
\centering
\subfigure[]{\epsfig{file=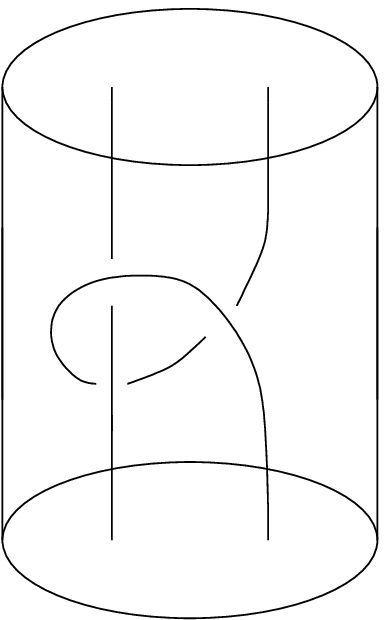, height=2cm}}
\hspace{1.5cm}
\subfigure[]{\epsfig{file=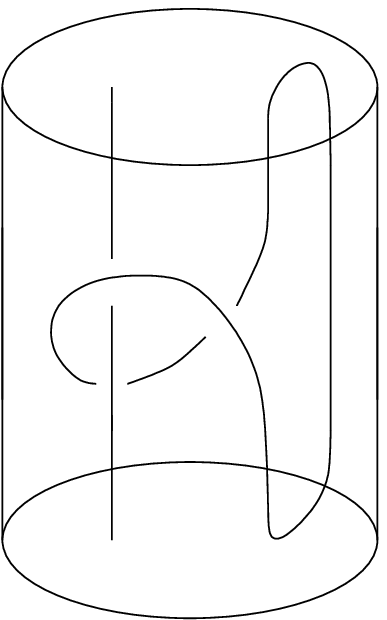, height=2cm}}
\caption{A string-link and its deformation closure.}
\label{stringlink}
\end{figure}

\smallskip
We now turn our attention to defining string-links in a homology ball $B^M$, these are the type of tangles we will be concerned with.

Fix a collection of points $p_1 < p_2 < p_3 < \cdots $ on the $x$-axis of $D^2$.  These induce  sets of points on $\partial ( \d2i )$ which we call the \textit{standard  points}.
By a \textit{string-link of n components}  $\sigma \subset B^M$  we mean a tangle
$\sigma : \cup_{i=1}^{n} I_i \rightarrow B^M $ such that under the identification of the boundary with $\partial (\d2i)$ we have  $\sigma |_{I_i} (j) = p_i \times j$, for $j=0,1$, where $p_i$ is the $i$-th standard point.

The two sets of {\em standard points} of a string link are the subsets of standard points on $D^2 \times \{ j \}$, $j=0,1$ which are the boundary points of a string-link component.

Given a string-link $\sigma \subset \d2i$, we can change a given set $A$ of (interval) components of $\sigma$ into  $S^1$ components by constructing non-intersecting paths on $\pa{\d2i}$  between the two endpoints of each  interval component, and pushing these paths and the endpoints of the components of $A$ slightly into the interior of $\d2i$.  This gives a tangle which we call  the  {\em deformation closure}  of $\sigma$ with respect to $A$.  An example is given in Figure~\ref{stringlink}, where the deformation closure is with respect to the right hand component.

We can {\em represent} any string-link $\sigma^{\prime} \subset B^M$ by a string-link $\sigma \subset \d2i$ with a specified set of framed surgery components, where $\sigma^{\prime}$ is obtained by carrying out surgery on the deformation closure of the surgery components of $\sigma$.

\smallskip

A {\em parenthesization}  of a set of standard points is a bracketing of that set (for example $(( p_1 p_2 )( p_3 ( p _4 p_5 )))$.  We call the parenthesization $(((( p_1 p_2 ) p_3 ) \cdots ) p_n )$ the {\em canonical parenthesization}.
A {\em parenthesization}  of a string-link  is a parenthesization of its two sets of standard points.

Note that a parenthesization on a string-link in $\d2i$ induces one on any string-link in $B^M$ it represents.

\begin{definition}
A {\em manifold string-link} is a  canonically  parenthesized, coloured, fram\-ed string-link in $\d2i$ with a set of linking components $X_L$ and surgery components $X_M$.
\end{definition}

We say that a manifold string-link is {\em regular} if the linking matrix of its surgery components is invertible (so  surgery yields a string-link in  a rational homology ball).

\smallskip

We will now define some actions on the set of string-links which we will make use of later.  As these are well known and somewhat fiddly, to define we gloss over the technical details and rely upon the readers intuition.

 If two string-links $\sigma_1$ and $\sigma_2$ in $\d2i$ have the same number of components  we may form a product $\sigma_1 \cdot \sigma_2$ in the usual way by ``putting $\sigma_2$ on top of $\sigma_1$''.
 If the string links are parenthesized or coloured we require that the parenthesization or colourings match on the two disks  identified under the composition.

We also define $\sigma_1 \otimes \sigma_2$ to be the string link obtained by ``placing $\sigma_2$ to the right of $\sigma_1$''.

Let $T \subset B^M$ be an $X$-coloured tangle and let $A \subset X$.
 Define $\varepsilon_A (T)$ to be the tangle obtained from $T$ by deleting all of the components with colours in $A$.
Further let $B$ be a set disjoint from $X$  and let $S \subset X \times B$.
We define $D_{S}(T)$ to be the coloured tangle obtained from $T$ by, for each $(x,b) \in S$, doubling the $x$-coloured component and colouring the double with $b$.
When dealing with string-links we may have to isotope them so that the endpoints lie on the appropriate standard points.

\section{Milnor's $\mu$-invariants}
Recall that given a ring $R$, a {\em R-homology sphere} is a 3-manifold $M$ such that $H_{q}(M;R) =H_{q}(S^3 ;R)$, for all integers $q$.  Similarly a {\em R-homology ball } is a 3-manifold $B^M$ with boundary $\partial (B^M ) = S^2$ such that $H_{q}(B^M;R) =H_{q}(B^3 ;R)$, for all  $q$, where $B^3$ is the 3-ball.
If  $R= \mathbb{Z}$  we do not specify the ring and just write {\em homology sphere} or {\em homology ball}.

Let $\sigma$ be a $l$-component string-link in a homology ball $\bsig$ with a fixed identification of $\partial (\bsig)$ with $\partial (\d2i )$ and let  $N(P)$ be a regular neighbourhood of the set of standard points and  $N( \sigma )$  a regular neighbourhood of the string-link. Then there are two inclusion maps
 $i_j :D^2 - N(P) \hookrightarrow \bsig  - N( \sigma )$ for $j=0,1$,
 where the map $i_j$ sends $D^2 - N(P)$ to the image of $D^2 \times \{ j \} - N(P)$ under the identification of  $\partial (\bsig)$ with $\partial (\d2i )$.
We use these $i_j$ to induce certain group isomorphisms as follows.

Let  $G$ be any group. The \textit{lower central series}, $G_q$ is defined inductively by $G=G_1$ and  $G_{q+1} = [G,G_q]$.

\smallskip

\noindent {\bf Stallings' Theorem (\cite{St:65}).}{\em
\: Let $\map{h}{A}{B}$ be a homomorphism of groups, inducing an isomorphism $\h{1}{A} \cong \h{1}{B} $ and an epimorphism from $\h{2}{A}$ onto $\h{2}{B}$.  Then, for finite $q \geq 0$, $h$ induces an isomorphism $ \quoq{A} \cong \quoq{B}$.
}

\smallskip

A Mayer-Vietoris calculation and a standard  application of Stallings' Theorem gives the following result.

\begin{proposition} \label{prop:isoms} {\em
$(i_j)_*$, $j=0,1$, induces isomorphisms
\[ \fquo{\fun{D^2 - N(P)}} \cong \fquo{\fun{\bsig - N(\sigma)}}. \]}
\end{proposition}

Let $F(l)$ be the free group on generators $x_1, \ldots , x_l$.  We will also denote the image of $x_i$ in the quotient group $\quoq{F(l)}$  by $x_i$ and the induced maps on the lower central series coming from Proposition~\ref{prop:isoms}  by $(i_j)_*$, $j=0,1$.
Since  we can identify $\fun{D^2 - N(P)}$ with $F(l)$, we have isomorphisms
\[ \fquo{F(l)} \overset{(i_0)_*}{\longrightarrow} \fquo{\fun{\bsig - \sigma}} \overset{(i_1)_*}{\longleftarrow} \fquo{F(l)}. \]
The composition $(i_1)_{*}^{-1} (i_0)_*$ gives a map
$ SL(l) \rightarrow Aut(\quoq{F(l)})$, where  $SL(l)$ is the set of string-links of $l$ components in a given homology ball $\bsig$.
It is not difficult to see that we in fact get a map
\smallskip
\begin{center}
$ \art_q : SL(l) \rightarrow Aut_0 (\quoq{F(l)}), $
\end{center}
\smallskip
where $Aut_0 (\quoq{F(l)})$ is the subgroup of $Aut (\quoq{F(l)})$ consisting of all automorphisms which map $x_i$ to a conjugate of itself and leaves the product $x_1 x_2 \cdots x_l$ of the generators fixed.
We call the map $\art_q$ the $q$-th \textit{Artin representation}.

The $i$-th \textit{longitude} $\lambda_i \in \quoq{F(l)}$ of a string-link $\sigma$ is defined in the following way.
Take a double of the $i$-th component of the string-link. This determines an element in the fundamental group of the complement.  Under $(i_1)_*^{-1}$ this gives an element in $\quoq{F(l)}$ which we call the longitude.
We have
\[ \art_q (\sigma)(x_i)=\lambda_i x_i \lambda_{i}^{-1}, \]
where $\lambda_i \in \quoq{F(l)}$ is the $i$-th longitude of $\sigma$.

Note that our longitudes are determined by the (black-board) framing and are not necessarily null-homologous.  It is easy to modify the content of this paper should we insist that the longitudes are null-homologous.

\begin{definition} We say that a string-link $\sigma$ has \textit{Milnor filtration n}, if all its longitudes are trivial in $\frac{F(l)}{F(l)_n}$.
\end{definition}

\smallskip

Let $\mathcal{P}$$(l)$ be the ring of formal power series in non-commuting variables $\ser{X}{l}$. The {\em Magnus expansion} is the homomorphism
\smallskip
\begin{center}
$ \mu : F(l) \rightarrow \mathcal{P}$$(l) $
\end{center}
\smallskip
defined on the generators of the free group by $\mu (x_i) = 1 + X_i$.

\begin{definition}
The {\em $\mu$-invariants}  of a string-link $\sigma$ in an integral homology ball are the coefficients of the monomials in the $X_i$ of the Magnus expansion of the $i$-th longitude $\lambda_i \in \fquo{F(l)}$.  Explicitly, the $i$-th $\mu$-invariant of {\em length} $n+1$  is
\[
\mu_{j_1 ,j_2 ,   \ldots , j_n ;i} = \text{Coeff}( X_{j_1} X_{j_2} \cdots X_{j_n} , \mu (\lambda_i))
\]
where $n \leq q$ and $\lambda_i \in \fquo{F(l)}$.
\end{definition}

It is well known that the longitudes $\lambda_i$ of $\sigma$ are trivial in $\frac{F(l)}{F(l)_n}$, that is $\lambda_i$ is of Milnor filtration $n$, if and only if all $\mu$-invariants of length $\leq n$ vanish.

\section{The Algebras} \label{sec:algebras}
The algebras we need are amalgamations of the usual algebras $\mathcal{A}$ and $\mathcal{B}$ from the theory of finite-type invariants (see~\cite{BN:95:2}).

\begin{definition}
Let $X,Y$ be finite disjoint sets.
 Then $\algyz$ is the space of formal $\mathbb{Q}$-linear combinations of uni-trivalent graphs whose trivalent vertices are oriented and whose univalent vertices are either coloured by elements of a set $Y$ or lie on the oriented coloured 1-manifold $(\cup_{x \in X} I_x )$,  which is called the \textit{skeleton},  modulo the STU, IHX and AS relations shown in Figure~\ref{fig:relations}.
\end{definition}

\begin{figure}
\[
\begin{array}{c}
\mathrm{STU:} \hspace{0.5cm} \STU{1cm}{0.5cm} \hspace{1cm} \mathrm{AS:} \hspace{0.5cm} \AS{1cm}{0.5cm} \\
\vspace{0.3cm} \\
\mathrm{IHX:} \hspace{0.5cm} \IHX{1cm}{0.5cm}
\end{array}
\]
\caption{The STU, AS and IHX relations.}
\label{fig:relations}
\end{figure}

Note that we allow trivalent graphs and the possibility that $Y=\emptyset$.

We denote the subspace of $\algyz$ where  every non-empty graph has a univalent vertex and all univalent vertices lie on the skeleton by $\mathcal{A}$$(\uparrow_X)$ and the subspace $\algez$ by $\mathcal{B}$$(Y)$.

The \textit{degree} of a uni-trivalent diagram is half of its number of vertices.

We say  an element of $\algyz$ is {\em connected} if it is a $\mathbb{Q}$-linear combination of connected uni-trivalent graphs.

\smallskip

\begin{figure}
\begin{center}
\epsfig{file=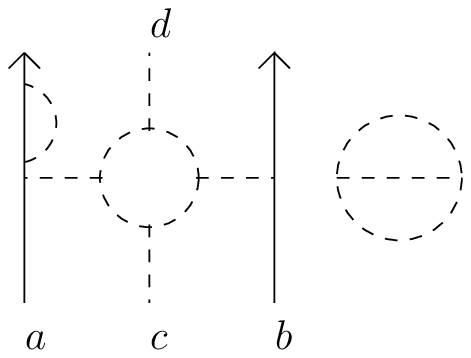, height=2cm}
\raisebox{10mm}{$\cdot$}
\epsfig{file=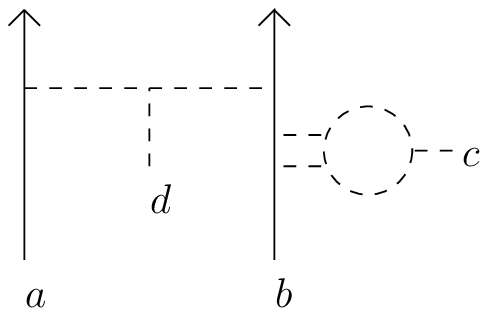, height=2cm}
\raisebox{10mm}{$=$}
\epsfig{file=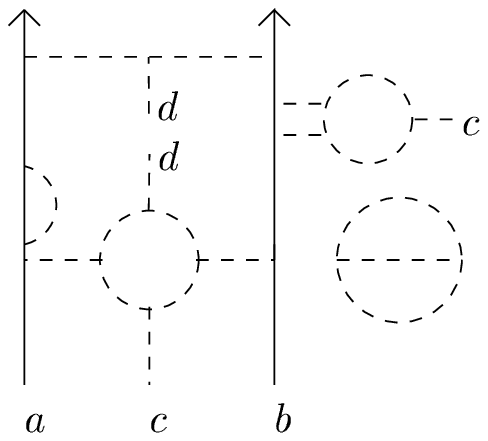, height=2cm}
\end{center}
\caption{An example of multiplication.}
\label{exmult}
\end{figure}

Let $D_1, D_2 \in \algyz$ then there is a  {\em product} $D_1 \cdot D_2$ given by the linear extension of the process of stacking the skeleton of $D_1$ on top of $D_2$ in such a way that the colours of the two skeletons match and taking the disjoint union of any trivalent components.
An example of the multiplication is given in Figure~\ref{exmult}.

There is also a notion of a coproduct $\cp$ in $\algyz$ which is the obvious extension of the usual coproduct of $\mathcal{A}$ (see~\cite{BN:95:2}). This makes  $\algyz$ into a graded co-commutative Hopf algebra where the grading is by the degree.
We denote the degree $n$ part by $\algyzn $ and, by abuse of notation, its graded completion again by $\algyz$.
 The primitives (ie. the elements such that $\cp (D) = 1 \otimes D + D \otimes 1$) of the Hopf algebra are the connected elements.

\smallskip

We will now look at some maps between these algebras. These properties hold since they hold in $\mathcal{A}$ and $\mathcal{B}$.

\begin{figure}
\[ \chiso_{\{b\}}
\left( \begin{array}{c}
 \epsfig{file=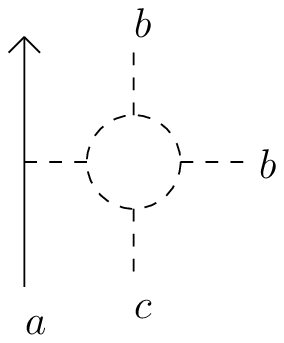, height=1.8cm}
\end{array} \right)
= \frac{1}{2}
\left( \begin{array}{ccc}
 \epsfig{file=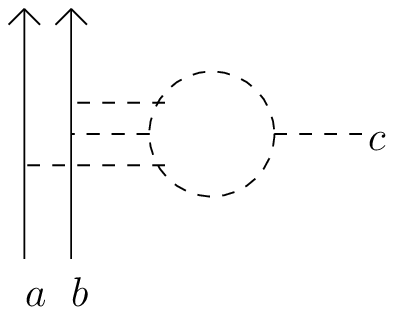, height=1.8cm}
 & \raisebox{0.9cm}{+}
  & \epsfig{file=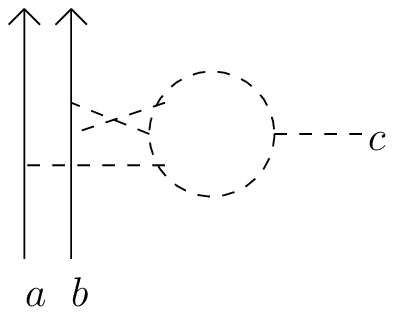, height=1.8cm}
\end{array} \right)
\]
 \caption{An example of the map $\chi_{ \{ b \} }$. }
\label{exchi}
\end{figure}

Let $Y^{\prime} \subset Y$. Define a map
\smallskip
\begin{center}
$\chiso_{ Y^{\prime} } : \algyz \rightarrow \algchiimage$
 \end{center}
\smallskip
by the linear extension of the process of adding $Y^{\prime}$ coloured skeleton components and
taking the average of all ways of placing the $ Y^{\prime}$ labeled univalent vertices on $\uparrow_{ Y^{\prime}}$.
See Figure~\ref{exchi} for an example of this map.

In fact $\chiso$ descends to a coalgebra isomorphism and we denote its inverse by $\sigma$.

If $X^{\prime} \subset X$, $Y^{\prime} \subset Y$  and $A=X^{\prime} \cup Y^{\prime}$. The map
$\varepsilon_{A} : \algyz \rightarrow \mathcal{A}$$( \uparrow_{X -X^{\prime}}, Y -Y^{\prime})$
is defined by setting every uni-trivalent graph with a uni-valent vertex on a $X^{\prime}$ coloured skeleton component or  with a $Y^{\prime}$ coloured vertex equal to zero.

\begin{figure}
\[ D_{\{(a,b) , (c,d) \}}
\left( \begin{array}{c}
 \epsfig{file=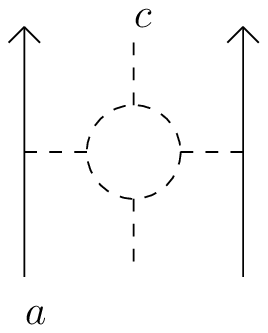, height=1.8cm}
\end{array} \right)
=
\begin{array}{c}
 \epsfig{file=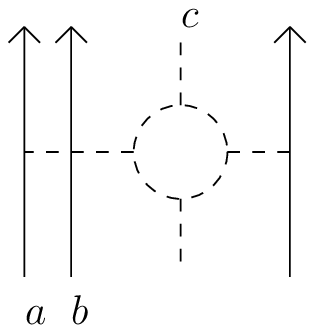, height=1.8cm}
 \raisebox{0.9cm}{+}
 \epsfig{file=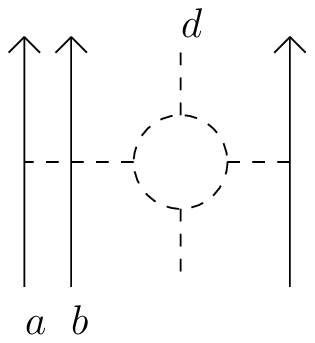, height=1.8cm}
 \raisebox{0.9cm}{+}
 \epsfig{file=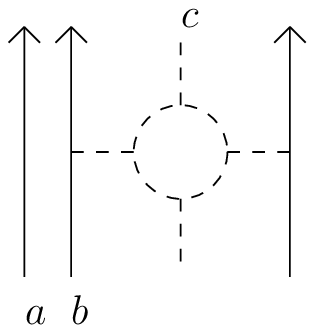, height=1.8cm}
 \raisebox{0.9cm}{+}
 \epsfig{file=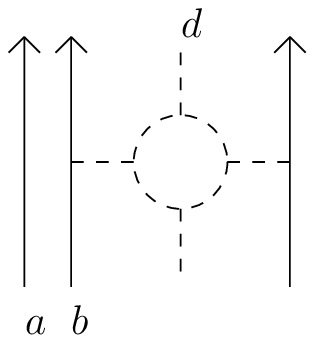, height=1.8cm}
\end{array}
\]
 \caption{An example of the map $D_{S}$.}
\label{exd}
\end{figure}

Let $B$ be some set  disjoint from both $X$ and $Y$ and let
$S \subset (X \cup Y) \times B$. Define
$D_{S} $ to be the linear extension of the operation which to each element $(a,b) \in  S$
 either, if $a$ is the label of a skeleton component, gives the sum of all ways of lifting the vertices lying on the $a$-coloured component  over    this component and its $b$-coloured double   and, if $a$ is the colour of a univalent vertex, is the sum of all ways of substituting the colour $a$ by $b$  (see figure~\ref{exd} for an example).

\medskip

We will be interested in two particular quotients of $\algyz$ which were  defined in \cite{ BN:95} and \cite{ HM:00} for the algebra $\mathcal{A}$.

Define $\algyzt$ to be the quotient of $\algyz$ by the ideal generated by all relations which set non-simply connected uni-trivalent graphs  equal to zero.  The connected elements are called  {\em trees} .
We will denote the connected part (ie. the primitives) of $\mathcal{B}$$^t (Y)$ by $\mathcal{C}$$^t(Y)$.

Also define $\algyzh$ to be the quotient of $\algyzt$ by the ideal generated by all relations which set connected uni-trivalent graphs with more than one univalent vertex either lying on the same skeleton component  or being labeled by the same colour, equal to zero.  We call this quotient the {\em homotopy quotient}.

It follows from \cite{BN:95}  that $\algyzh$ is a quotient of $\algyzt$ and $\chiso$ descends to  isomorphisms.

\smallskip

There is a well known map (see eg. \cite{HM:00, GL:P1}) which relates trees to Lie algebras.
Let $ \text{Lie} (l) = \oplus_{n \geq 1} \text{Lie}_n (l)$, be the free $\mathbb{Q}$ Lie algebra on $l$ generators $\ser{X}{l}$.
Also let $\mathcal{C}$$^t(Y,a)$ be the subspace of $\mathcal{C}$$^t(Y \cup \{ a \})$ consisting of connected elements in which every uni-trivalent graph has exactly one univalent vertex coloured by some $a \notin Y$.

\begin{figure}
\[
\epsfig{file=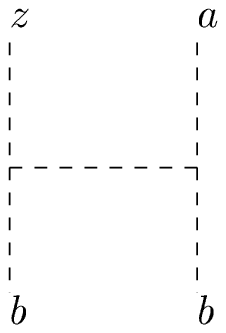, height=2cm}
\raisebox{0.9cm}{$\rightarrow$}
\epsfig{file=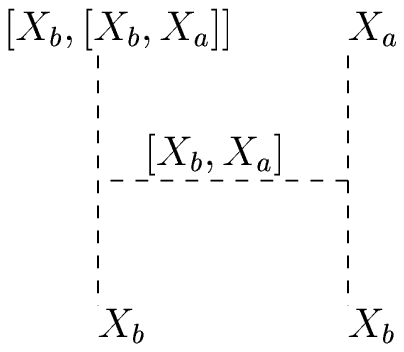, height=2cm}
\raisebox{0.9cm}{$\rightarrow [X_b ,[X_b , X_a]]$}
\]
\caption{An example of the isomorphism
$\mathcal{C}$$^t_3 (\{ a,b \},z) \rightarrow \text{Lie}_3 (2)$.}
\label{fig:primlie}
\end{figure}

Fix a bijection between the colouring set $Y$ and the generators $\ser{X}{l}$ of the free Lie algebra, where $|Y|=l$. Then given some element $D \in \mathcal{C}$$^t_n(Y,a)$  label the edges ending in a $Y$-coloured univalent vertex with the corresponding generator of the Lie algebra.
Now assign an element of the Lie algebra to each unlabelled edge according to the rule that whenever an unlabelled edge meets two  edges labelled by $X$ and $X^{\prime}$ in $ \text{Lie} (l)$ (in the direction of the orientation) assign the commutator $[X,X^{\prime}]$ to that edge.  This labels the edge coloured by $a$ and we take this to be our element of $\text{Lie}_n (l)$.
See Figure~\ref{fig:primlie} for an example.
It is not hard to see that this  gives an isomorphism from $\mathcal{C}$$^t_n(Y,a)$ to $\text{Lie}_n (l)$.

\begin{figure}
\[
\epsfig{file=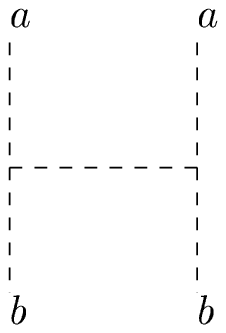, height=2cm}
\raisebox{0.9cm}{$\rightarrow$}
\hspace{2mm}
\epsfig{file=mufig/is1, height=2cm}
\raisebox{0.9cm}{+}
\hspace{2mm}
\epsfig{file=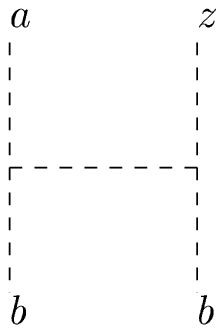, height=2cm}
\raisebox{0.9cm}{$\rightarrow$}
\hspace{2mm}
\raisebox{0.9cm}{2}
\hspace{1mm}
\epsfig{file=mufig/is1, height=2cm}
\raisebox{0.9cm}{$\rightarrow 2[X_b ,[X_b , X_a]]$}
\]
\caption{An example of
$j_z : \mathcal{C}$$^t_3 (\{ a,b \}) \rightarrow \text{Lie}_3 (2)$.}
 \label{fig:mapj}
\end{figure}

Finally, we define a map
$ j_y : \mathcal{C}$$^t_n(Y) \rightarrow \text{Lie}_n (l)$ for $y \in Y$
by summing over all of the ways replacing exactly one of the $y$-coloured vertices by some $a \notin Y$ and then using the above map to get an element in $\text{Lie}_n (l)$.
An example is given in Figure~\ref{fig:mapj}.

\bigskip

Although we are assuming some familiarity with the Kontsevich integral, we briefly review some relevant properties here.

\begin{itemize}
\item For our purposes the Kontsevich integral $Z$ is an $\mathcal{A}$$(\uparrow_{X})$ valued universal finite-type invariant of $X$-coloured framed parenthesized tangles and the degree $n$ part of $Z$, $Z_n$  is a degree $n$ finite-type invariant.

\item Let $\pi^{h} : \mathcal{A}$$(\uparrow_{X}) \rightarrow \mathcal{A}$$^h(\uparrow_{X})$ be projection. Then by \cite{BN:95}, $\pi^h \circ Z$ is a well defined invariant of ($X$-coloured framed parenthesized) string-links  up to link-homotopy, where {\em link-homotopy} is an equivalence relation which allows ambient isotopy and each component of the tangle to pass through itself.

\item Let $T$ and $T^{\prime}$ be tangles then $Z(D_{S}(T)) = D_{S}(Z(T))$,
$Z(\varepsilon_{A}(T)) = \varepsilon_{A}(Z(T))$ and $Z(T \cdot T^{\prime}) = Z(T) \cdot Z(T^{\prime})$.

\item $Z(T) \in \mathcal{A}$$(\uparrow_{X})$ is group-like and so can be written as $\exp ( C )$ where $C \in \mathcal{C}$$(X)$ is connected.
\end{itemize}

\section{The \AA rhus Integral} \label{sec:aarhus}

The \AA rhus integral, $Z^M$,  was introduced by Bar-Natan, Garoufalidis, Rozansky and Thurston in the series of papers \cite{BGRT:AI, BGRT:AII, BGRT:AIII} as a universal finite type invariant of rational homology 3-spheres.
In this series it was remarked that it  extends to an invariant of  links in rational homology spheres.
In this section we define \AA rhus integral.
The reader is referred to the \AA rhus trilogy  for a thorough exposition  of the invariant.

\medskip

The {\em pre-normalized \AA rhus integral} of regular manifold string-links  $\PZM$ is defined by the following composition:
\smallskip
\begin{center}
$\PZM : RMSL \overset{\z}{\longrightarrow} \Alm \overset{\sigma_{X_M}}{\longrightarrow} \alglm \overset{\fgi}{\longrightarrow} \algle$
\end{center}
\smallskip
where:
\begin{itemize}
\item  $RMSL$ is the set of regular manifolds string-links with linking components coloured by $X_L$ and  surgery components coloured by $X_M$.
\item $\z \overset{def}{=} \nu^{\otimes |X_L \cup X_M |} \cdot D_{\{-\}\times X_L \cup X_M} (\nu ) \cdot Z$, is the Kontsevich integral as normalized in \cite{LMMO:99}.
\item $\fgi$ is  {\em formal Gaussian integration} with respect to the variables $X_M$. It is described below.
\end{itemize}

\begin{definition}
The {\em \AA rhus integral} of a regular manifold string-link, $\sml$ is given by
\[
\ZM(\sml ) = \PZM (U_+ )^{- \sigma_{+}} \cdot \PZM (U_- )^{- \sigma_{-}} \cdot  \PZM(\sml )
\]
where $\sigma_{\pm}$ is the number of $\pm$ve eigenvalues of the linking matrix of $\varepsilon_{X_L} (\sml) $ and $U_{\pm}$ is the unknot with framing $\pm 1$.
\end{definition}

We will now go on to define formal Gaussian integration.
Let $D_1, D_2 \in \algrs$, define
\[
\bracy{D_1}{D_2} =
\left(
\begin{array}{l}
\text{sum of all ways of gluing all legs labeled } \\ y \text{ on } D_1 \text{ to the legs labeled } y \text{ on }  D_2, \\ \text{for all the colours } y \in Y .
\end{array}
\right) ,
\]
 where this sum is non-zero only if the number of $y$-coloured legs of $D_1$ equals the number of $y$-coloured legs of $D_2$, for all $y \in Y$.

\begin{figure}
\[
\left\langle \begin{array}{c}
 \raisebox{0.7cm}{$\frac{1}{2}$}
 \hspace{1mm}
 \epsfig{file=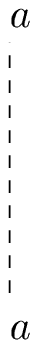, height=1.4cm}
\raisebox{0.7cm}{$+$}
\hspace{2mm}
\epsfig{file=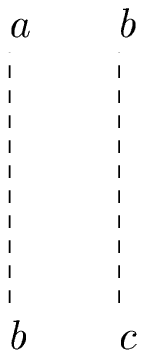, height=1.4cm}
\hspace{2mm}
,
\hspace{2mm}
 \epsfig{file=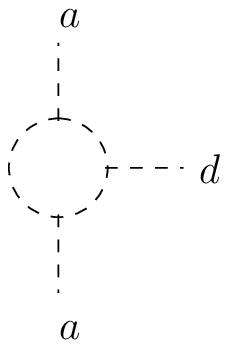, height=1.4cm}
 \raisebox{0.7cm}{$+$}
 \epsfig{file=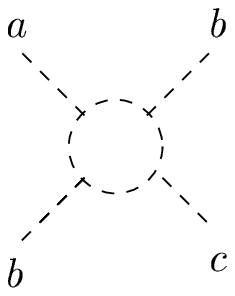, height=1.4cm}
\end{array} \right\rangle_{ \{ a,b,c \} }
\begin{array}{c}
\raisebox{0.7cm}{$=$}
 \hspace{3mm}
 \epsfig{file=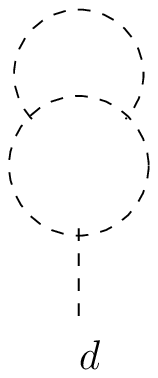, height=1.4cm}
 \hspace{1mm}
 \raisebox{0.7cm}{$+$}
 \hspace{2mm}
 \raisebox{0.7cm}{$2$}
 \hspace{1mm}
 \epsfig{file=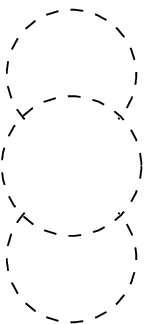, height=1.4cm}
\end{array}
\begin{array}{c}
\raisebox{0.7cm}{$=$}
 \hspace{3mm}
 \raisebox{0.7cm}{$2$}
 \hspace{1mm}
 \epsfig{file=mufig/b5.eps, height=1.4cm}
\end{array}
\]
 \caption{An example of $\bracs{D_1}{D_2}$.}
\end{figure}

\medskip

It is a well known and easily seen fact that for an $X$-coloured tangle $T$, the Kontsevich integral $\z $ may be written in the form
\[ \sigma (\z (T)) =
\exp ( \sum_{x,y \in X} \frac{1}{2}l_{xy} \; x \frown y ) + (\text{other stuff}),\]
where $(l_{xy})$ is the linking matrix of $T$.  Recall that the degree one elements (which look like $\frown$) are called {\em struts}.

Therefore, given a regular manifold string-link $\sml$, with linking components $X_L$ and surgery components $X_M$, we can write
\[
\sigma_{X_M} \z (\sml) = \exp (\sum_{x,y \in X_M} \frac{1}{2} l_{xy} \; x \frown y) \cdot P \overset{def}{=} \exp (Q/2)  \cdot P,
\]
where now $(l_{xy})$ is the linking matrix of $\varepsilon_{X_L} (\sml) $.
Since $\sml$ is regular, $(l_{xy})$ is invertible and so we can define:
\[
Q^{-1} = \sum_{x,y \in X_M} l^{xy} \; x \frown y
\]
where $(l^{xy})$ is the inverse matrix of $(l_{xy})$.

Writing $\sigma_{X_M} \z (\sml) \overset{def}{=}  \exp (Q/2)  \cdot P$, we define {\em formal Gaussian integration} as:
\[
\fgi  \exp (Q/2)  \cdot P = \bracxm{\exp (-Q^{-1} /2)}{P}.
\]

\medskip

It  is known (\cite{BGRT:AII}) that $\PZM$ is invariant under isotopy and a handle slide of any component around a surgery component, and  $\ZM$ is invariant under stabilization on the surgery components and so $\ZM$ descends to an invariant of string-links in a rational homology ball.
Summarizing this we have:

\begin{proposition} {\em
The \AA rhus integral $Z^M$ is an invariant of framed parenthesized string-links in rational homology balls.}
\end{proposition}

\medskip

At this point we fix some notation.
Let $\pi^t : \algxy \rightarrow \algxyt$ be projection.  Then
$\ZMT \overset{def}{=} \pi^t \circ \ZM$,
$\bracy{-}{-}^t \overset{def}{=} \pi^t \circ \bracy{-}{-}$
and so on.
We use similar notation for the projection $ \pi^h : \algxy \rightarrow \algxyh$.

\section{The \AA rhus Integral and the $\mu$-invariants}

Let $\sml$ be a manifold string-link with the canonical parenthesization such that  the determinant of the linking matrix of the surgery components is $\pm 1$ (so $\sml$ represents a string-link in an integral homology ball).
Further, for convenience, we set $X_L = \{1, \ldots , l \}$ ,  $X_M = \{l+1, \ldots , l+m \}$ and assume that the components of the manifold string-link have numerically increasing colours from left to right.
We call such a manifold string-link {\em nice}.

At times we will need to add an extra  linking component to the manifold string-link.  We will add this component to the left of the others and colour it with $0$.
We denote the new colouring set $X_L \cup \{ 0 \}$ by $X_L +1$.

The extra $0$-coloured component is going to correspond to a longitude of the string-link and as such is only considered up to link homotopy.  Consequently, rather than working with the algebra $\algtl1e$, we add an additional homotopy relation on the colour $0$, and we call the resulting algebra $\alghtl1e$.

Given a set of colours $X$, let  $1_X = \otimes_{x \in X} 1$, be the trivial tangle coloured by $X$. When $X$ contains only one element, $x$ say,  we will just write $1_x$.

In this section we consider the longitudes as elements
$\lambda_i = \lim_{ \hspace{-0.6cm}\raisebox{-1mm}{$\longleftarrow$}} \lambda_i^{(n)}$ of the nilpotent completion
$ \widehat{F(l)} = \lim_{ \hspace{-0.6cm}\raisebox{-1mm}{$\longleftarrow$}} \quo{F(l)}$,
where $\lambda_i^{(n)} \in \quo{F(l)}$.

The reader is referred to \cite{HM:00} for the motivation behind the formula in the following proposition.

\begin{proposition} \label{mu:lambda} {\em
Let $\sml$ be a nice  manifold string-link   and let  $\lambda_i$, $1 \leq i \leq l$, be its $i-th$ longitude regarded as a pure braid. Then
\begin{equation} \label{mu:form}
Z^{M;h,t}(\lambda_i \otimes 1_{X_M}) = \pi^{h,t} (Z^M (1_0 \otimes \sml )^{-1} (D_i Z^M ( \sml ))^{b_i})
\end{equation}
where $b_i = Z^M (\beta_i)$ and  $\beta_i$ is the braid coloured by $\{0, \ldots , l+m \}$ inducing the permutation $(i-1 \; i-2 \cdots  1 \; 0)$, $a^b$ denotes the conjugation  $bab^{-1}$, $D_i = D_{\{(i,0)\}}$ and $\pi^{h,t}$ is projection onto $\alghtl1e$.}
\end{proposition}

\begin{Remark}
In formula~\ref{mu:form} we are assuming that $\sml$, $\lambda_i \otimes 1_{X_M}$ and $1_0 \otimes \sml$ have the canonical parenthesization and $\beta_i$ has the canonical parenthesization on the bottom and the `$i$-th double of the canonical parenthesization' on the top.
\end{Remark}

\begin{Remark}
 Since $\Zmht_0 = \Zmht$ we need only consider the pre-normalized \AA rhus integral. Also note that $\z^t =Z^t$.
\end{Remark}

\medskip

We need a few technical lemmas to prove the proposition.

\begin{lemma} \label{mu:double} {\em
Let $\sml$ be a manifold string-link, $i \in X_L$ and $D_i = D_{\{(i,0)\}}$. Then
\[Z^{M}(D_i (\sml )) = D_i(Z^{M}(\sml )).\] }
\end{lemma}

\begin{proof}
\[
\begin{split}
\pzm{D_i ( \sml )} &=  \sideset{}{_{X_M}^{FG}}\int  \sigma_{X_M} \z (D_i ( \sml)) \\
                &=  \sideset{}{_{X_M}^{FG}}\int  D_i( \sigma_{X_M}( \z( \sml ))) \\
                &=  D_i \left( \sideset{}{_{X_M}^{FG}}\int  \sigma_{X_M} \z (\sml ) \right) \\
                &= D_i(\pzm{ \sml})
\end{split}
\]
where the second equality is a standard property of the Kontsevich integral.
The third follows since $i \in X_L$ and the formal Gaussian integration is with respect to the variables $X_M$.

The result follows since $D_i$ respects multiplication.
\end{proof}

\begin{lemma} \label{mu:glue} {\em
Let $P \in \algyz$ be of degree $n$ and contain no struts both of whose univalent vertices are coloured by elements of $Y$, and let $Q \in \algyz$ consist entirely of struts coloured by $Y$.  Then if $\bra{Q}{P}{Y}^t$ is non-zero, it is of degree at least $ n-[ \frac{n}{2} ]$, where  $[a]$ denotes the integer part of $a$.}
\end{lemma}

\begin{proof}
In order to determine the minimum possible degree of $\bra{Q}{P}{Y}^t$ it is enough to find a graph  $R \in \mathcal{A}$$_n (\uparrow_X , Y)$ which satisfies the conditions  the lemma imposes on $Q$ and will allow the maximum number of struts to be glued in under $\bra{-}{R}{Y}^t$, whilst still resulting in a non-zero term.
It is not hard to see that such an example  occurs when $R$ consists entirely of  struts with exactly one $Y$-colored vertex and, if $n$ is odd, one strut both of whose vertices are $X$-coloured.
In such a situation we can glue in at most $[ \frac{n}{2} ]$ struts. The result follows upon noting that each strut we glue in reduces the degree by one.
\end{proof}

Note that if $P$ and $Q$ are as in the above lemma and $\bra{Q}{P}{Y}$ is non-zero, then it is of degree at least $ n-[ \frac{3n}{4} ]$.

\begin{figure}
\begin{center}
 \epsfig{file=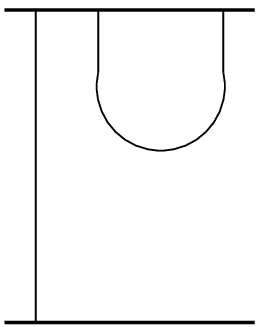, height=1.8cm}
\hspace{1cm}
\epsfig{file=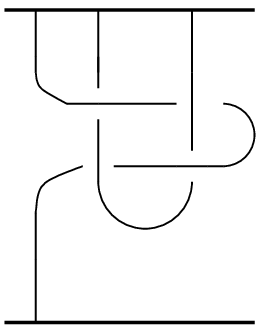, height=1.8cm}
\caption{Two tangles which differ by a pure braid.}
\label{mu:differ}
\end{center}
\end{figure}

\begin{definition}
We say that two tangles  $T$ and $T'$ {\em differ by a pure braid $p \in PB_{n+1}$} if $T'$ can be obtained from $T$ by replacing a copy of $D^2 \times I$ which intersects $T$ in a trivial string-link with the pure braid $p$  (see Figure~\ref{mu:differ}).
\end{definition}

We need the following result of Stanford.
\begin{theorem}[Stanford \cite{St:96}]{\em
Let  $T$ and $T'$ be two tangle which differ by a pure braid $p \in PB_{n+1}$.  Then for any finite type invariant, $v$, of degree less than $n+1$ we have $v(T)=v(T')$.}
\end{theorem}

\begin{lemma} \label{mu:uni} {\em
Suppose Y is the disjoint union of compact 1-manifolds $Y_0 , Y_L , Y_M$ where $Y_M$ consists entirely of copies of $S^1$.  Further suppose that $T_i$ ,  $i=1,2$ , are two tangles which agree on $Y_L \cup Y_M$ and on each component of $Y_0$ the maps differ by an element in the lower central series $ \pi_1 (M-T_i |_{Y_L} )_{n+1}$ , where $M$ is the manifold obtained by surgery on $T_i |_{Y_M}$.  Then the images of $\pzm{(T_i)}$ in $\algspecial$, where the homotopy filtration is on the $Y_0$ components, agree. }
\end{lemma}

\begin{proof}
First note that $\pi_1 (\d2i - T_i |_{Y_L \cup Y_M}  )$ is generated by the meridians of $T_i |_{Y_L \cup Y_M}$.
This means that  there exists a ball which intersects $T_i |_{Y_L \cup Y_M}$ in a trivial string-link such that the map of the fundamental groups induced by the inclusion of the trivial string-link into $T_i |_{Y_L \cup Y_M}$ is surjective.

Now since $ \pi_1 (M-T_i |_{Y_L} )$ is also generated by the meridians of $T_i |_{Y_L \cup Y_M}$ and $T_1$ and $T_2$ differ by elements in $ \pi_1 (M-T_i |_{Y_L} )_{n+1}$, we see that $T_2$ can be obtained from $T_1$ by handle sliding around the $Y_M$ components, modifying the $Y_0$ components by homotopy and modifying $T_1$ (inside the ball described above) by pure braids in
$F(|Y_L \cup Y_M|)_{n+1} \subset F(|Y_L \cup Y_M|) \subset \text{PB}(|Y_L \cup Y_M|+1) $.

Stanford's theorem tells us  that the modification by the pure braids does not affect finite-type invariants of degree less than $n+1$ and since the homotopy relations are applied to the $Y_0$ components so tangles differing under these two moves have the same image under $\z_{\leq n}^{h,t}$.
Finally, formal Gaussian integration takes care of the handle slides and the result then follows  by Lemma~\ref{mu:glue}.
\end{proof}

\noindent {\bf Proof of Proposition~\ref{mu:lambda}.}
Let $\lambda_{i}^{(n)}$ be a representative of the longitude $\lambda_i$ in $F(l)/F(l)_{n+1}$ which we regard as a pure braid of $l+1$ components.
Now the two tangles
$(1_0 \otimes \sml)(\lambda_{i}^{(n)} \otimes 1_M)$ and $\beta_i D_i(\sml) \beta_i^{-1}$
both represent the union of the manifold string-link, $\sml$, and the longitude and therefore, by Stallings' Theorem, satisfy the conditions of Lemma~\ref{mu:uni}.
Then
\[
Z^{M;h,t}_{< n -[\frac{n}{2} ]}((\lambda_{i}^{(n)} \otimes 1_M) \cdot (1_0 \otimes \sml ) )=
Z^{M;h,t}_{< n -[\frac{n}{2} ]}(\beta_i \cdot D_i (\sml ) \cdot \beta_{i}^{-1}).
\]
Now since $\lambda_i^{(n)} \in F(l) / F(l)_{n+1}$, it is a word in the meridians of the linking components only, so $\lambda_{i}^{(n)}$ only has crossings with the linking components of $\sml$. Therefore
 \[Z^{M;h,t}_{< n -[\frac{n}{2} ]}((\lambda_{i}^{(n)} \otimes 1_M) \cdot (1_0 \otimes \sml ) ) = Z^{M;h,t}_{< n -[\frac{n}{2} ]}(\lambda_{i}^{(n)} \otimes 1_M) Z^{M;h,t}_{< n -[\frac{n}{2} ]} (1_0 \otimes \sml ).\]
Similarly,
 \[Z^{M;h,t}_{< n -[\frac{n}{2} ]}(\beta_i \cdot D_i (\sml ) \cdot \beta_{i}^{-1}) =
 Z^{M;h,t}_{< n -[\frac{n}{2} ]}(\beta_i )Z^{M;h,t}_{< n -[\frac{n}{2} ]}(D_i (\sml ) )Z^{M;h,t}_{< n -[\frac{n}{2} ]}(\beta_{i}^{-1} ).\]

Finally solving for  $Z^{M;h,t}_{< n -[\frac{n}{2} ]}((\lambda_{i}^{(n)} \otimes 1_M)$, and letting $n$ tend to infinity gives the result.\; \; \qed

\medskip

Having found a formula for the \AA rhus integral of the longitudes we turn our attention to finding a formula for the Magnus expansion of the longitudes.

\begin{definition}
An {\em expansion} is a homomorphism  $J:F(l) \rightarrow \mathcal{P}$$(l)$ such that
$J(x_i) = 1 +X_i + (\text{higher order terms})$
, where $F(l)$ is the free group on the generators $x_1, \ldots , x_l$ and  $\mathcal{P}$$(l)$ is the ring of formal power series in  non-commuting variables $X_1, \ldots , X_l$.
\end{definition}

Clearly the Magnus expansion is an expansion in this sense.
We show that the left hand side of formula~\ref{mu:form} can be regarded as an expansion and then we apply the following result of Lin to write the Magnus expansion of the longitudes in terms of the \AA rhus integral.

\begin{lemma}[Lin \cite{Li:97}] \label{linslemma} {\em
Let $J$ be any expansion and $\mu$ be the Magnus expansion.  Then there exist a unique unipotent automorphism $\Psi : \mathcal{P}$$(l)\rightarrow \mathcal{P}$$(l)$ such that $\mu = \Psi \circ J$. }
\end{lemma}

Recall that a map $\Psi$ is said to be {\em unipotent} if for all $a \in  \mathcal{P}$$(l)$ of degree $n$, $\Psi(a)=a + O(n+1)$.

\medskip

 $\Ahtlo$ is a graded co-commutative Hopf algebra whose space of primitives is isomorphic to $\Cht{X_L +1}$, the space  of connected elements of $\Bhtlo$.
 Let $X_i \in \Ahtlo$ denote the element of degree one which has a single chord between the skeleton components coloured by 0 and $i$.

Then $X_1 , \ldots , X_l$  generate a free non-commutative power series ring
$\mathcal{P}$$(l) = \mathcal{P}$$(X_1 , \ldots , X_l) \subset \Ahtlo, $
(since the primitives of $\Ahtlo$ are isomorphic to
$\Cht{X_L +1}$ which is naturally decomposed as
 $ \Ct{X_L,0} \oplus \Ct{X_L} = \text{Lie}(l) \oplus \Ct{X_L}$
and the first summand corresponds to $\mathcal{P}$$(X_L)$ by the isomorphism described in Section~4).

\begin{figure}
\begin{center}
 \epsfig{file=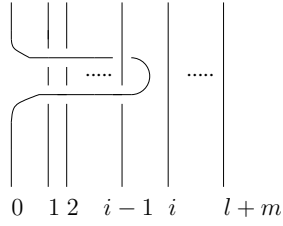, height=3cm}
\end{center}
\caption{The pure braid $ \sigma_{0,i} \otimes 1_M$.}  \label{mu:inj}
\end{figure}

Let  $SL(X_M )$ denote the monoid of string-links in $\d2i$ which are coloured by $X_M$ and $PB(l+1)$ be the pure braid group on $l+1$ generators.
There is map $\iota: F(l) \rightarrow PB(l +1) \otimes SL(X_M)$ defined by the formula
$x_i \mapsto \sigma_{0,i} \otimes 1_M$
where $\sigma_{0,i}$ is the generator of the pure braid group which wraps the  $0$-th strand once around the $i$-th as in Figure~\ref{mu:inj}.
The composition of this  with $ \Zmht $ gives a map $J:F(l) \rightarrow \alghtl1e$.

\begin{lemma} {\em
The map $J :F(l) \rightarrow \alghtl1e  $ defined above is an expansion.}
\end{lemma}

\begin{proof}
Let $x_i \otimes 1_M$ denote the generators of $F(l) \otimes id \subset PB(l+1) \otimes SL(X_M)$.
Then
\[
\zmht{x_i} = i \circ \zht{p_1(x_i)} = \Exp{Y_i}
\]
where $Y_i = X_i + O(2) \in \Ahtl1m$,
$i:\Ahtlo \rightarrow \Ahtl1m$ is  inclusion and $p_1: PB(l+1) \otimes SL(X_M) \rightarrow PB(l+1)$ is  projection onto the first component.

To prove the lemma we have to show that the image of $J$ lies in $\mathcal{P}$$(l)$.
But this follows since every diagram in $\zmht{x}$, where $x \in F(l) \otimes id$, must have a vertex lying on the skeleton component coloured $0$ since the removal of the $0$-coloured component trivializes the braid.
Thus $Y_i \in \Ct{X_L ,0} = \text{Lie}(l)$.
\end{proof}

Now applying Lemma~\ref{linslemma} to Proposition~\ref{mu:lambda} gives:
\begin{theorem} \label{mu:mu} {\em
Let $\sml$ be a nice manifold string-link.  Then
\[
\mu (\lambda_i) = \Psi \circ \pi^{h,t} (Z^M (1_0 \otimes \sml )^{-1} (D_i Z^M ( \sml ))^{b_i})
\]
where $b_i = Z^M (\beta_i)$ and  $\beta_i$ is the braid coloured by $\{0, \ldots , l+m \}$ inducing the permutation $(i-1 \; i-2 \cdots  1 \; 0)$, $a^b$ denotes the conjugation  $bab^{-1}$, $D_i = D_{\{(i,0)\}}$, $\pi^{h,t}$ is projection onto $\alghtl1e$ and $\Psi$ is a unipotent automorphism.}
\end{theorem}

\begin{Remark}
Habegger and Masbaum's theorem in \cite{HM:00} relating the $\mu$-invariants of string-links in $\d2i$ to the Kontsevich integral (which is obviously contained in Theorem~\ref{mu:mu}) holds in the algebra $\Ahoxl1$, which we  define to be $\Axl1$ with the homotopy relation applied to the colour 0, so there are no tree relations on the $X_L$-coloured components (this is true since \cite{HM:00}'s~Lemma~12.5 only requires the homotopy relation and in their Lemma~12.6, the homotopy relation ensures that the appropriate elements are trees).
However it is interesting to note that attempts by the author to remove the ``$t$'' in this section failed as the normalized Kontsevich integral $\z$ does not respect multiplication, necessitating the descent into $\alghtl1e$.
\end{Remark}

\section{The First Non-vanishing $\mu$-invariant}
As an application of the above we generalize  the Habegger-Masbaum formula by expressing  the first non-vanishing Milnor invariants in terms of the first non-vanishing term of the tree part of the \AA rhus integral. We note that Habegger gave a different proof of this result in \cite{Ha2000}.

\medskip

Recall that $ \text{Lie} (l) $ is the free $\mathbb{Q}$ Lie algebra on $l$ generators $\ser{X}{l}$.
There is a canonical graded isomorphism of
$\oplus_{n\geq 1}(F(l)_n / F(l)_{n+1} ) \otimes \mathbb{Q}$ with $ \text{Lie} (l)= \oplus_{n \geq 1} \text{Lie}_n (l)$.
Now if $\sigma$ has Milnor filtration $n$, we can consider the longitudes $\lambda_i$ as elements in $F(l)_n / F(l)_{n+1}$ and we denote the corresponding element in $\text{Lie}_n (l)$ by $\MU{i}{n}$.
We call the $\MU{i}{n}$  the {\em Milnor invariants of degree n}.

\begin{theorem} \label{mu:formulae} {\em
Let $\sml$ be a nice manifold string-link representing a string-link $\sigma$. Then

(i) $\Zmht (\sml) = 1+ O(n)$ if and only if the string-link represented by $\sml$ is of Milnor filtration $n$,

(ii) the first non-vanishing Milnor invariants of the string-link $\sigma$ determine and are determined by the first non-vanishing term of $\zmt{\sml} -1$ through the Habegger-Mausbaum formula:
\[ \mu_{i}^{(n)} (\sml) = j_i (\xi )  \]
where $\zmt{\sml} = 1 + \xi + O(n+1)$ and $j_i : \mathcal{C}$$^{t}_{n} (X_L) \rightarrow \text{Lie}_n (X_L) $ is the map described in Section~4. }
\end{theorem}

\begin{proof}

\noindent (i)
First suppose that $\zmt{\sml} = 1+O(n)$, then
$\pht (D_i(\zm{\sml})) = 1+O(n) $ giving $\pht ((D_i(\zm{\sml}))^{b_i}) = 1+O(n)$
(since the lower degree terms of the conjugating $b_i$'s cancel).

Also we have  $\pht (\zm{1_0 \otimes \sml}^{-1}) = 1 +O(n)$.

Since multiplication can not reduce the degree and $\Psi$ is unipotent we get
\[
\mu (\lambda_i) = \Psi \circ \pi^{h,t} (Z^M (1_0 \otimes \sml )^{-1} (D_i(Z^M ( \sml )))^{b_i})  = 1+ O(n),
\]
and the result follows since all $\mu$-invariants of length $\leq$ n vanish if and only if $\lambda_i$ is trivial in $\pi(\bsig - \sigma) / \pi(\bsig - \sigma)_n$, where $\sigma$ is a tangle represented by $\sml$.

\smallskip

Conversely, suppose that $\lambda_i$ is trivial in $\pi(\bsig - \sigma) / \pi(\bsig - \sigma)_n$.
 Then
\[\coeff{x_{\iota_1} x_{\iota_2} \cdots x_{\iota_{n-r-1}} }{\mu (\lambda_i )} = 0  \text{, \: for } 0 \leq r \leq n-1,\]
 and so
\[
\mu (\lambda_i) = \Psi \circ \pi^{h,t} (Z^M (1_0 \otimes \sml )^{-1} (D_i(Z^M ( \sml )))^{b_i})  = 1+ O(n).
\]
As $\Psi$  is unipotent it follows that
\[
\pi^{h,t} (Z^M (1_0 \otimes \sml )^{-1} (D_i(Z^M ( \sml )))^{b_i})  = 1+ O(n).
\]
Thus $D_i(\zmht{\sml}) = 1 + O(n)$ and so $\zmt{\sml} = 1 + O(n)$.

\medskip

\noindent (ii)
Suppose that the first non-vanishing $\mu$-invariant is of degree $n$. Then by the above
$\zmht{\sml} = 1 +\xi +O(n+1)$, where $\xi$ is of degree $n$.
By Proposition~\ref{mu:mu}, Lemma~\ref{mu:double} and the unipotency of  $\Psi$  we have
\[
\mu (\lambda_i) = \Zmht_{\leq n}(1_0 \otimes \sml )^{-1}
\cdot \Zmht_{\leq n} (\beta_i)
\cdot \Zmht_{\leq n} ( D_i(\sml ))
\cdot \Zmht_{\leq n} (\beta_i)^{-1} + O(n+1).
\]
Since $\zm{\beta_i}$ can be written as the exponential  of a sum of connected elements, this can be written as
\[
\begin{split}
&=(1 - 1_0 \otimes \xi)(1+\zeta)(1+ \pi^{h,t} D_i (\xi))(1-\zeta) + O(n+1) \\
&= 1+ \pi^{h,t} D_i (\xi)+ O(n+1)
\end{split}
\]
where $\zeta = \Zmht_n (\beta_i)$ and $\xi$ is as in the statement of the theorem.

Looking at the degree $n$ part of this formula we see that the terms of $ D_i (\xi)$ which do not have a vertex on the 0-coloured skeleton component cancel with the terms of $1_0 \otimes \xi$, and any terms of $ D_i (\xi)$ with more than one vertex on the 0-coloured skeleton component are killed off by the projection $\pi^{h,t}$.
So what remains is an element of $\mathcal{C}$$^t_n(X_L, 0) = \text{Lie}_n (l)$ and it is easy to see that this is exactly the element $j_i (\xi )$.

Finally, the determined by part follows since $j_i$ is injective (see \cite{HM:00}).
\end{proof}

\noindent{\bf Acknowledgements}
This work was carried out while the author was visiting Georgia Institute of Technology, who he would like to thank for their hospitality.  He would also like to thank his supervisor Stavros Garoufalidis for suggesting this problem.

\end{document}